\documentclass{amsart}         
\usepackage{amssymb}         
\usepackage{latexsym}         
\newcommand{\ZZ}{{\mathbb Z}}         
\newcommand{\RR}{{\mathbb R}}         
\newcommand{\CC}{{\mathbb C}}         
\newcommand{\NN}{{\mathbb N}}

\newtheorem{theorem}{Theorem}            
\newtheorem{remark}{Remark}            
\newtheorem{lemma}{Lemma}[section]            
\newtheorem{prop}[lemma]{Proposition}            
\newtheorem{coro}[lemma]{Corollary}            
\newtheorem{definition}[lemma]{Definition}            
\renewcommand\qedsymbol{$\Box$}         
         
\newcounter{smalllist}

\newcommand{\rin}{r_{\rm in}}     
\newcommand{\rout}{r_{\rm out}}     
\newcommand{\Pset}{{\mathcal P}}     
\newcommand{\Qdot}{\dot{Q}}     
\newcommand{\Pdot}{\dot{P}}     
\newcommand{\Pbox}{{\mathcal{P}}_{\rm b}}     
\newcommand{\Pir}{{\mathcal{P}}(r)}     
\newcommand{\Pinfty}{{\mathcal{P}}(\infty)}     
\newcommand{\Fbar}{\overline{F}}     
\newcommand{\dist}{{\rm dist}}     
\newcommand{\worte}{\mathcal{W}}     
\newcommand{\boxen}{\mathcal{B}}     
     
\begin{document}         
\title[Linear repetitivity]{Linear repetitivity, I. Uniform subadditive   
  ergodic theorems and  applications}         
\author[D.~Damanik, D.~Lenz]{David Damanik$\,^{1,2}$, Daniel Lenz$\,^{2}$}         
\maketitle         
\vspace{0.3cm}             
\noindent             
$^1$ Department of Mathematics 253--37, California Institute of Technology,             
Pasadena, CA 91125, U.S.A.\\[2mm]             
$^2$ Fachbereich Mathematik, Johann Wolfgang Goethe-Universit\"at,             
60054 Frankfurt, Germany\\[2mm]       
E-mail: \mbox{damanik@its.caltech.edu, dlenz@math.uni-frankfurt.de}\\[3mm]             
2000  AMS Subject Classification: 52C23, 37A30, 37B50     
     
\begin{abstract}     
This paper is concerned with the concept of linear repetitivity in the   
theory of tilings. We prove a general uniform subadditive ergodic theorem for linearly repetitive tilings. This theorem unifies and extends various known (sub)additive ergodic theorems on tilings. The results of this paper can be applied in the study of both random operators and lattice gas models on tilings.      
\end{abstract}         
         
\section{Introduction}         
In a recent paper, Lagarias and Pleasants studied linearly and densely     
repetitive tilings \cite{LP}. It was shown that these structures are diffractive and they proposed to consider linearly repetitive tilings as models of ``perfectly ordered quasicrystals.''   
   
In fact, several special classes of  linearly repetitive tilings have attracted much attention. One such class is given by tilings arising from primitive substitutions. They have been studied in several contexts \cite{GH,GS,hof,LuP,sol1,sol2}, including random Schr\"odinger operators and lattice gas models. Both the study of lattice gas models and the study of random Schr\"odinger operators require a uniform subadditive ergodic theorem. The appropriate theorem has been established in \cite{GH}. In the one-dimensional case there is another important class of examples of linearly repetitive structures, namely, Sturmian dynamical systems whose rotation number has bounded continued fraction expansion. Again, this class allows for a uniform subadditive ergodic theorem. This has been shown by one of the authors \cite{len} (cf.~\cite{dl2, l3} for applications). These results immediately raise the following question:   
   
\begin{itemize}     
\item[(Q)] Does linear repetitivity imply a uniform subadditive ergodic theorem?     
\end{itemize}     
This question is answered in the affirmative by Theorem     
\ref{subadditive} in Section \ref{Uniform} of this paper (cf. \cite{l3}   
as well). This theorem   generalizes the theorem of   
\cite{len}. Moreover,  combined   
with the known linear repetitivity of tilings generated by primitive   
substitution \cite{DZ,DHS,sol2}, it gives a conceptual proof for the subadditive ergodic theorem of \cite{GH}. Of course, this theorem also implies an additive ergodic theorem. However, this additive ergodic theorem is not as effective as the corresponding theorem of \cite{LP}, as it does not contain an error estimate (cf. Section \ref{Uniform}).   
   
We defer discussion of the methods used in the proofs of our results to the   
corresponding sections. However, we would like to emphasize the   
following perspective in our considerations: Our point of view is a   
purely local one. Thus, the key object of our studies is neither a   
tiling nor a species of tilings but rather certain sets of pattern   
classes. The appropriate sets are defined in Definition \ref{admissible}   
and termed admissible. The advantage of this point of view is   
twofold. Firstly, in this approach, the uniformity of results is built in as the local structure is uniform for all tilings in the species. Secondly, the role of asymptotic translation invariance appearing in the subadditive ergodic theorems is clarified (cf. Section \ref{Specializing}).    
     
The article is organized as follows. In Section \ref{Preliminaries} we      
review basic facts on tilings and fix some notation. Section   
\ref{Uniform} contains a rather general form of a subadditive ergodic   
theorem. This is the main result of this paper. It gives an affirmative   
answer to Question (Q). In Section \ref{Specializing} we specialize   
the main theorem to various situations. This recovers several known   
(sub)additive ergodic theorems. Finally, in Section \ref{Applications}, we sketch applications of the foregoing results in the study of random operators associated to tilings.

\section{Preliminaries}\label{Preliminaries}     
The aim of this section is to introduce certain notions and to fix some     
notation.      
     
Consider a set consisting of subsets of $\RR^d$ which are homeomorphic to the closed unit ball in $\RR^d$ and pairwise disjoint up to their boundaries. Such a set of sets will be called a pattern if it is finite. It will be called a tiling (of $\RR^d$) if the union of its elements equals the     
whole space. Its elements will be called tiles. For certain applications, it is useful to consider decorated tiles and patterns. A pattern with decorations from a set $\Gamma$ is a set $M$ of pairs $m=(a_m,c_m)$ with $a_m\subset \RR^d$ and $c_m\in \Gamma$ such that $\{a_m:m\in M\}$ is a pattern. One should think of a pair $(a_m,c_m)$ as a tile colored or decorated by $c_m$. The following definitions apply to both patterns and decorated patterns. However, to avoid tedious repetitions, they are phrased in terms of patterns only. If a pattern $M$ is contained in a pattern or tiling $N$, we write $M\subset N$ and say that $M$ is a subpattern of $N$. Similarly, if a tile $t$ belongs to a pattern or tiling $N$, we write $t\in N$. For a pattern $M$, we define the underlying set $s(M)$ by      
$$     
s(M)=  \cup_{m\in M} m\subset \RR^d.     
$$      
The inner radius $\rin(M)$ of a pattern $M$ is defined by     
$$     
\rin(M)= \max\{r\in \RR: \exists x\in \RR^d, K(x,r)\subset     
s(M)\},     
$$     
and the outer radius $\rout(M)$ of a pattern $M$ is defined by     
$$     
\rout(M)=\min\{r\in \RR: \exists x\in \RR^d, K(x,r)\supset     
s(M)\},     
$$     
where $K(x,r)$ denotes the closed ball around $x$ with radius $r$. The     
existence of the minimum and maximum in question follows by compactness     
of $s(M)$. For a pattern $M$ and a closed set  $B$ homeomorphic to the   
unit ball  with $s(M) \supset B$, define the restriction $M\cap B$ of $M$   
to $B$ by    
   
\begin{equation}\label{boxenschnitt}   
M\cap B=\{ m\cap B: m\in M, m\cap {\rm int}(B)\neq \emptyset\}.   
\end{equation}   
In the applications we have in mind, $B$ will be either a box   
(cf.~Section \ref{Uniform}) or a closed ball.   
   
There exists a natural equivalence relation on the set of patterns. Two patterns are equivalent if and only if they agree up to translation. The class of a pattern will also be called a pattern class or an abstract pattern. Similarly, an abstract tile is the class of a tile up to translation. The relations ``$\in$'' and ``$\subset$'' (resp., the functions $\rin$ and $\rout$) give rise to relations (resp., functions) on abstract patterns in the obvious way. The induced relations (resp., functions) will be denoted by the same symbols. Similarly, concepts such as connectedness of patterns, disjointness, or distance of tiles in patterns, etc.~can easily be carried over to abstract patterns. This will tacitly be done in the sequel, whenever necessary. Moreover, we will sometimes omit the word abstract in abstract patterns if no confusion can arise.    
     
Our point of view is a purely local one. Thus, the following definition introduces the main object of our studies.      
     
\begin{definition}\label{admissible}     
A set $\Pset$ of abstract patterns in $\RR^d$ is called admissible if it satisfies the following conditions.     
     
\begin{itemize}     
\item[(i)] $P\in \Pset$, $Q\subset P$ implies $Q\in \Pset$.     
\item[(ii)] There exist $0< r_{\rm min},r_{\rm max} < \infty$ with $r_{\rm min}\leq \rin(a)\leq \rout(a)\leq r_{\rm max}$ for all abstract tiles $a\in \Pset$.     
\item[(iii)] Let $P\in \Pset$ with representative $\Pdot$ with $0\in \Pdot$ and $r>0$ be given. Then, there exists a $Q\in \Pset$ with representative $\Qdot$ with  $K(0,r) \subset s(\Qdot)$ and $\Pdot \subset \Qdot$.      
\end{itemize}     
\end{definition}     
     
In the sequel we will be exclusively concerned with admissible sets $\Pset$. For each admissible set, there is  a natural set of tilings associated with it. Conversely, to a tiling of $\RR^d$, one can associate a set of abstract patterns. This is the content of the next definition.     
     
\begin{definition}\label{associated}     
{\rm (a)} Let $T$ be a tiling of $\RR^d$. The set $\Pset(T)$ of abstract patterns associated to $T$ is defined to be the set of classes of subpatterns of $T$.\\     
{\rm (b)} Let $\Pset$ be an admissible set of abstract patterns. A tiling     
$T$ is said to be  associated to $\Pset$ if $\Pset(T)\subset \Pset$.\\   
{\rm (c)} Let $\Pset$ be an admissible set of abstract patterns. The set    
of all tilings $T$ associated to $\Pset$ with the topology induced by   
the metric   
$$   
{\rm d}(T,S)= \inf \left\{\epsilon: T \cap B(0,\tfrac{1}{\epsilon}) = (S+t)\cap B(0,\tfrac{1}{\epsilon}), t\in \RR^d, \|t\|\leq \epsilon \right\}   
$$   
is a topological space denoted by $\Omega(\Pset)$ {\rm (}cf.~\cite{sol1}{\rm )}.   
\end{definition}     
     
This article is centered around the notion of linear repetitivity. This     
notion has been studied in \cite{LP} for Delone sets in $\RR^d$. In our     
context it is given in the following definition.     
     
\begin{definition}\label{linearrepetitive}     
An admissible set $\Pset$ of patterns is called linearly repetitive if     
there exists a constant $c_{{\rm LR}}>0$ such that every $P\in \Pset$ with     
$\rout(P)\geq 1$ is contained in every $Q\in \Pset$ with $\rin(Q)\geq     
c_{{\rm LR}} \cdot \rout(P)$.      
\end{definition}     
   
An important property of the tiling space associated to a linearly repetitive $\Pset$ is the following:   
   
\begin{prop}\label{compact}   
If the admissible $\Pset$ is linearly repetitive, then $\Omega(\Pset)$ is compact.   
\end{prop}   
\noindent{\it Proof.} This follows by rather standard arguments once it is realized that linear repetitivity implies finiteness of the number of pattern classes with a prescribed maximal outer radius. For the reader's convenience, we include a proof in Appendix A. \hfill \qedsymbol   
   
\medskip   
     
Let us finish this section by discussing the role of Delone sets and the      
Voronoi construction in our context. Recall that a subset $D$ of $\RR^d$ is called a Delone set if there exist positive constants $r_{0}$ and $r_{1}$ such that each ball in $\RR^d$ of radius at least $r_1$ contains a point of $D$ and each ball of radius at most $r_0$ does not contain more than one point of $D$. The Voronoi construction assigns to each $x$ in a given Delone set $D$ the set $V(x)=\{y\in \RR^d: \dist(x,y)\leq \dist(z,y), z\in D\}$, where $\dist(\cdot,\cdot)$ denotes Euclidean distance. Then, $V(D)=\{V(x):x\in D\}$ is a tiling of $\RR^d$ by convex polytopes (cf.~\cite{sen}). Proposition 5.2 of \cite{sen} says that $\Pset(V(D))$ is admissible for a Delone set $D$. Thus, Delone sets give rise to admissible sets. This motivates the following definition.     
     
\begin{definition}     
A Delone set $D$ is called linearly repetitive if $\Pset(V(D))$ is linearly repetitive.      
\end{definition}     
     
\begin{remark}   
{\rm Using the following proposition, it is not hard to show that this definition of linear repetitivity for Delone sets agrees with the definition of \cite{LP}.}      
\end{remark}     
     
\begin{prop}\label{Delone} Let $D$ be a Delone set. Then for each $x\in D$, the tile $V(x)$ is determined by the points of $D$ lying inside a ball of     
radius $2 r_1$ around $x$.      
\end{prop}     
\noindent{\it Proof.} This is just Corollary 5.1 in \cite{sen} \hfill     
\qedsymbol

\section{The Main Theorem}\label{Uniform}     
This section is devoted to a proof of a rather general uniform subadditive ergodic theorem. The proof is similar to that of \cite{len}, which in turn uses ideas of \cite{GH} (cf.~\cite{dl2,l3} for further details). The formulation relies on patterns on boxes. Thus, we will start this section with a discussion of boxes.      
     
A box $B$ in $\RR^d$ is a subset of the form $B=\{(x_1,\ldots, x_d):     
a_j\leq x_j\leq b_j, j=1,\ldots,d\}$, where $a_j<b_j\in \RR$ for each     
$j$. The length of the $j$-th side is denoted by $l_j$, that is,     
$l_j=b_j-a_j$. The volume and the surface area of a box $B$ are denoted by $|B|$ and  $\sigma(B)$, respectively. Moreover, let the width $\omega(B)$ of a box $B$ be defined by $\omega(B)=\min\{l_j:j=1,\ldots,d\}$. For $r\in \RR^+$, an $r$-box is a box whose sidelengths satisfy     
$$     
r\leq l_j\leq 2 r, \,\; j=1,\ldots, d.     
$$     
The set of all boxes (resp., $r$-boxes) is denoted by $\boxen(\RR^d)$   
(resp., $\boxen(r)$). A box-pattern (resp., $r$-pattern) is a pattern   
$M$, where $s(M)$ is a box (resp., $r$-box). For a box $B$ and a pattern   
(or tiling) $M$ with $s(M)\supset B$, the box-pattern derived from $M$   
by restricting to $B$ denoted by $M\cap B$ has been defined in   
(\ref{boxenschnitt}).      
     
Now, let an admissible set of abstract patterns $\Pset$ be given. The set     
$\Pbox$ of abstract box-patterns derived from $P$ consists of all abstract patterns $Q$ which have representatives $\Qdot$ of the form $\Qdot=\Pdot\cap B$, where $B$ is a box and $\Pdot$ is a representative of $P\in \Pset$. If     
$B$ is an  $r$-box, the abstract pattern $Q$ is called an abstract     
$r$-pattern. The set of all abstract $r$-patterns derived from $\Pset$ is     
denoted by $\Pir$. Moreover, let $\Pinfty$ be defined by      
$$     
\Pinfty= \bigcup_{r>0} \Pset(r).     
$$      
The functions $l_j$, $|\cdot|$, $\sigma$, and $\omega$ induce functions on $\Pbox$ in the obvious way, which will be denoted by the same symbols. The inclusion relation $\subset$ on the set of boxes induces a relation on $\Pbox$, again denoted by $\subset$. That is, the relation $P\subset Q$ for $P,Q \in \Pbox$ holds if and only if there exist boxes $B_P,B_Q$ and representatives $\Pdot,\Qdot$ of $P$ and $Q$, respectively, with $s(\Pdot)=B_P$, $s(\Qdot)=B_Q$ and $\Pdot= \Qdot\cap B_P$.  Similarly, the equation     
$$     
P=\bigoplus_{j=1}^n P_j     
$$     
for $P,P_j\in \Pbox$, $j=1,\ldots,n$ is defined to hold if and only if there exist representatives $\Pdot$ of $P$ and $\Pdot_j$ of $P_j$, $j=1,\ldots,n$, with     
$$     
s(\Pdot)=\bigoplus_{j=1}^n s(\Pdot_j).     
$$     
Here, the equation $B=\bigoplus_{j=1}^n B_j$ for boxes $B,B_j$, $j=1,\ldots,n$ is defined to hold if and only if the $B_j$ are     
pairwise disjoint up to their boundaries and their union is $B$. Equations of the form $ P=\bigoplus_{j=1}^n P_j$ (resp., $B=\bigoplus_{j=1}^n B_j$) are called decompositions or partitions of patterns (resp., boxes). The notion of linear repetitivity appropriate to box-patterns is contained in part (ii) of the next proposition.     
     
\begin{prop}\label{LRbox}     
Let $\Pset$ be admissible and $\Pbox$ be as above. Then the following are equivalent:\\     
{\rm (i)} $\Pset$ is linearly repetitive, that is, there exists a constant $c_{\rm LR}$ such that every $Q\in \Pset$ with $\rout(Q)\geq 1$ is contained in every $P\in \Pset$ with $\rin(P)\geq c_{{\rm LR}} \cdot \rout(Q)$.\\     
{\rm (ii)} There exists a constant $c_{\rm LR,b}$ such that every $P\in \Pset(r)$ with $r\geq 1$ is contained in every $Q \in \Pset(c_{{\rm LR, b}} \cdot r)$.      
\end{prop}     
\noindent{\it Proof.} This is straightforward. \hfill \qedsymbol     
     
\medskip     
     
We can now introduce the class of subadditive functions.     
     
\begin{definition} Let $\Pset$ be admissible.\\     
{\rm (a)} A function $F:\Pbox \longrightarrow \RR$ is   
called subadditive if there exist nonnegative constants  $d_{\rm F}$ and   
$r_{\rm F}$ and a nonincreasing function $c_{\rm F}:[ r_{\rm F},   
\infty)\longrightarrow \RR$  with  $\lim_{r\to \infty} c_{\rm F}(r)=0$  such that     
\begin{itemize}     
\item[(i)] $F(P)\leq \sum_{j=1}^n F(P_j) + \sum_{j=1}^n c_{\rm   
    F}(\omega(P_j)) |P_j|$ for $P=\bigoplus_{j=1}^n P_j$ with $\omega(P_j)\geq r_{\rm F}$,      
\item[(ii)] $|F(P)|\leq d_{\rm F} |P|$.     
\end{itemize}     
{\rm (b)} A function $F:\Pbox\longrightarrow \RR\cup\{\pm \infty\}$ is called additive if  both $F$ and $-F$ are subadditive.      
\end{definition}     
   
We will be interested in means of subadditive functions. Our main result states the existence of a certain limit of means of this kind. These   
means are introduced in the next definition.   
   
\begin{definition}   
Let $F$ be subadditive on $\Pset$. For $r\geq r_F$ the means $F^+(r)$ and $F^-(r)$ are defined by   
$$   
F^+(r)= \sup\left\{ \frac{F(P)}{|P|}: P\in \Pset(r)\right\},\;\:\;\: F^-(r)= \inf\left\{ \frac{F(P)}{|P|}: P\in \Pset(r)\right\}.   
$$   
\end{definition}   
   
The following proposition is well known. In the context of subadditive     
functions on Delone sets, it was proved in \cite{LP}. For the     
convenience of the reader, we include a sketch of the proof.     
     
\begin{prop}\label{grenzwertfplus}     
Let $F$ be a subadditive function. Then, the following equation holds,     
$$     
\lim_{r\to \infty} F^+(r)=\inf_{r\geq r_{\rm F}}\{F^+(r)+      
c_{\rm F}(r)\}.     
$$      
\end{prop}     
\noindent{\it Proof.} Denote the infimum by $\Fbar$. We show (i) $\Fbar\leq \liminf_{r\to \infty} F^+(r)$ and (ii) $\limsup_{r\to \infty} F^+(r)\leq \Fbar$.     
     
(i) This is clear by      
$$     
\liminf_{r\to \infty} F^+(r)=\liminf_{r\to \infty}\left( F^+(r) +     
c_{\rm F}(r) \right)\geq  \inf_{r\geq r_{\rm F}}\{F^+(r)+      
c_{\rm F}(r)\}.     
$$     
     
(ii) Fix an arbitrary $r_0\geq r_{\rm F}$. Now, every $P\in \Pset(r)$     
with $r\geq r_0$ arbitrary can be written as a sum of patterns in     
$\Pset(r_0)$. The subadditivity condition together with a short     
calculation then implies     
     
\begin{equation}\label{stern}     
\frac{F(P)}{|P|}\leq F^+(r_0)+  c_{\rm F}(r_0).     
\end{equation}     
As $P\in \Pset(r)$ was arbitrary, equation (\ref{stern}) implies     
$$     
F^+(r)\leq F^+(r_0)+ c_{\rm F}(r_0)     
$$     
for all $r\geq r_0$. This proves (ii) and finishes the proof of the     
proposition.  \hfill \qedsymbol     
     
\medskip     
     
We can now prove the main theorem of this section.     
     
\begin{theorem}\label{subadditive}     
Let $\Pset$ be admissible and linearly repetitive and let $F$ be subadditive on $\Pbox$. Then the limits $\lim_{r\to \infty} F^+(r)$ and $\lim_{r\to     
\infty} F^-(r)$ exist and are equal. In particular, the equation     
$$     
\lim_{r\to \infty} F^+(r)=\lim_{|P|\to \infty, P\in \Pinfty}\frac{F(P)}{|P|}     
$$     
is valid.     
\end{theorem}      
\noindent{\it Proof.} This is proved by contraposition. So, assume     
$\liminf_{r\to \infty} F^-(r)< \limsup_{r\to \infty}F^+(r)$. Thus, by Proposition \ref{grenzwertfplus} there exist a $\delta>0$, a sequence     
$n(k)$ with $n(k)\to \infty$ for $k\to \infty$, and $Q_k\in \Pset(n(k))$     
with     
     
\begin{equation}\label{widerspruch}     
\frac{F(Q_k)}{|Q_k|}\leq F^+(n(k)) -\delta.     
\end{equation}     
W.l.o.g.~we can assume $n(k)\geq r_{\rm F}$. Choose an arbitrary $k\in     
\NN$ and consider some arbitrary $P\in \Pset(3 c_{\rm LR,b} n(k))$. Here      
$c_{\rm LR,b}$ is as defined in Proposition \ref{LRbox}. Let $\Pdot$ be     
an arbitrary representative of $P$ with underlying box $B=     
s(\Pdot)$. By partitioning each side of $B$ into three parts of equal     
length, the box $B$ can be decomposed into $3^d$ congruent smaller     
boxes, all belonging to $\boxen( c_{\rm LR,b} n(k))$. There is only one of these smaller boxes which does not intersect the boundary of $B$. Call it  $B_{\rm int}$. The decomposition of $B$ into smaller boxes induces a decomposition of $ \Pdot$ into $(c_{\rm LR,b} n(k))$-patterns. Denote by $\Pdot_{\rm int}$ the pattern with $s(\Pdot_{\rm int}) =B_{\rm int}$. By linear repetitivity, $\Pdot_{\rm int}$ contains a representative $\Qdot_k$ of $Q_k$. As the distance of $B_{\rm int}$ to the boundary of $B$ is bigger than or equal to $c_{\rm LR,b} n(k)$, the same is true for the distance of $s(\Qdot_k)$ to the boundary of $B$. Thus, $B$ can be written as     
$$     
B=\bigoplus_{j=0}^n B_j     
$$     
with suitable $B_j\in \boxen(n(k))$, $j=1,\ldots,n$, and $B_0=s(\Qdot_k)$.  This induces a decomposition of $P$ of the form $P=\bigoplus_{j=0}^n P_j$ with $P_j\in \Pset(n(k))$, $j=1,\ldots,n$, and $P_0= Q_k$. By subadditivity of $F$ and $(\ref{widerspruch})$ this implies     
     
\begin{eqnarray*}     
\frac{F(P)}{|P|}&\leq& \sum_{j=1}^n \frac{F(P_j)}{|P_j|}\frac{|P_j|}{|P|}     
+ \frac{F(Q_k)}{|Q_k|}\frac{|Q_k|}{|P|} +      
\sum_{j=0}^n c_{\rm F}(\omega(P_j)) \frac{|P_j|}{|P|}\\     
&\leq& \sum_{j=0}^n F^+(n(k))\frac{|P_j|}{|P|} - \delta  \,   
\frac{|Q_k|}{|P|} +  c_{\rm F}(n(k))  \\   
&\leq& F^+(n(k)) - \delta \, \frac{|Q_k|}{|P|} + c_{\rm F}(n(k)).     
\end{eqnarray*}     
Here we used the bound $\frac{F(P_j)}{|P_j|}\leq F^+(n(k))$, valid for arbitrary $P_j\in \Pset(n(k))$. Since $Q_k$ belongs to $\Pset(n(k))$ and $P$ belongs to $\Pset(3 c_{\rm LR,b} n(k))$, we can estimate     
$$     
\frac{|Q_k|}{|P|}\geq \frac{(n(k))^d}{(2 \cdot 3 c_{\rm LR,b} n(k))^d}=\frac{1}{(6  c_{\rm LR,b})^d}.     
$$     
Putting all this together, we arrive at      
$$     
\frac{F(P)}{|P|}\leq F^+(n(k)) - \frac{1}{(6  c_{\rm LR,b})^d} \, \delta +  c_{\rm F}({n(k)}).     
$$     
Since $P\in \Pset(3 c_{\rm LR,b} n(k))$ was arbitrary, this implies     
$$     
F^+(3 c_{\rm LR,b} n(k))\leq F^+(n(k)) - \frac{1}{(6  c_{\rm LR,b})^d} \, \delta +  c_{\rm F}(n(k)).     
$$     
As this holds for arbitrary $k\in \NN$, we can now take the limit on     
both sides using Proposition \ref{grenzwertfplus} and obtain $\Fbar\leq \Fbar -\delta \frac{1}{(6 c_{\rm LR,b})^d}$, a contradiction. This finishes the proof. \hfill \qedsymbol     
   
\medskip     
     
As a corollary we get an additive ergodic theorem. This is our version of Theorem 4.1 of \cite{LP}. Note, however, that we are not able to estimate the convergence rate (cf.~Remark \ref{Konvergenzrate} below).     
     
\begin{coro}\label{additive}     
Let $\Pset$ be linearly repetitive and let $F$ be an additive function on $\Pbox$. Then the following equation holds,     
$$     
\lim_{r\to \infty} F^+(r)=\lim_{\omega(P)\to \infty, P\in \Pbox}     
\frac{F(P)}{|P|}.     
$$      
\end{coro}     
\noindent{\it Proof.}  The decomposition technique of the proof of Proposition \ref{grenzwertfplus} applied to the subadditive function $F$  gives     
     
\begin{equation}\label{A}     
\limsup_{\omega(P)\to \infty, P\in \Pbox} \frac{F(P)}{|P|} \leq  \liminf_{r\to \infty}  F^+(r).     
\end{equation}     
Since $-F$ is subadditive as well, this equation immediately implies     
     
\begin{equation}      
\limsup_{\omega(P)\to \infty, P\in \Pbox} \frac{-F(P)}{|P|} \leq \liminf_{r\to \infty} (-F)^+(r).     
\end{equation}     
Multiplying by $(-1)$ and using $(-F)^+(r)=-F^-(r)$, we get     
     
\begin{equation} \label{B}     
\liminf_{\omega(P)\to \infty, P\in \Pbox} \frac{F(P)}{|P|} \geq \limsup_{r\to \infty} F^-(r).     
\end{equation}     
By $(\ref{A})$, $(\ref{B})$, and the foregoing theorem, the corollary     
follows. \hfill \qedsymbol     
     
\begin{remark}\label{Konvergenzrate}     
{\rm It is not possible to derive an estimate on the rate of convergence in the subadditive theorem. This can be seen from the following example. Let $f:\RR^+\longrightarrow [0,\infty]$ be an arbitrary monotonically decreasing function. Let $\Pbox$ be an arbitrary set of box-patterns derived from an admissible $\Pset$. Define $F:\Pbox\longrightarrow \RR$ by $F(P)= |P|     
f(|P|)$. As $f$ is decreasing, the function $F$ is subadditive. Moreover, we have $\frac{F(P)}{|P|}=f(|P|)$. Since $f$ was an arbitrary decreasing function, this shows that the rate of convergence in the subadditive ergodic theorem cannot be estimated.}     
\end{remark}     
   
\begin{remark}   
{\rm It appears that uniform subadditive ergodic theorems are special features of linearly repetitive structures. To support this, in the appendix we will exhibit examples of strictly ergodic structures for which a uniform subadditive ergodic theorem does not hold. The examples will be given by Sturmian subshifts whose rotation number has rapidly increasing continued fraction coefficients.}   
\end{remark}     
   
\section{Specializing the Main Theorem}\label{Specializing}   
In this section we derive various corollaries from the subadditive   
ergodic theorem. First, we consider (sub)additive functions on   
boxes on a conrete tiling or Delone set. We then use our methods to give   
a direct proof of the (known) unique ergodicity of dynamical systems   
arising from linearly repetitive tilings. Finally, we discuss how the   
theorems of \cite{GH,LP,len} fit into our context.   
   
Let us first introduce the appropriate notion of subadditivity and   
translation invariance.   
\begin{definition} Let  $w$ be a function on the set $B(\RR^d)$ of all   
  boxes in $\RR^d$.\\   
{\rm (a)} The function $w$ is called subadditive if there exists   
a constant $r_w$ and a nonincreasing function   
$c_w:[r_w,\infty)\longrightarrow \RR$ with $\lim_{r\to \infty} c_w(r)=0$  such that   
   
\begin{itemize}     
\item[(i)] $w(B)\leq \sum_{j=1}^n w(B_j) + \sum_{j=1}^n c_w(   
  \omega(B_j)) |B_j|$      
for $B=\bigoplus_{j=1}^n B_j$, with $\omega(B_j)\geq r_w$,     
\item[(ii)] $|w(B)|\leq d_w |B|$.     
\end{itemize}      
{\rm (b)} Let $T$ be a Delone set or a tiling. The function $w$ is called   
asymptotically  $T$-invariant if there exists a constant $r_w$ and a   
nonincreasing function $e_w:[r_w,\infty)\longrightarrow \RR$ with $ \lim_{r\to \infty}    
e_w(r)=0$ such that   
   
\begin{itemize}     
\item[(iii)] $|w(B)-w(B+t)|\leq  e_w(\omega(B))|B|$ if $(B \cap T) +t = (B+t)\cap T$ and $\omega(B)\geq r_w$.   
\end{itemize}     
\end{definition}   
   
We can now easily derive subadditive theorems for functions on Delone   
sets or tilings.   
   
\begin{coro}\label{tilingsubadditive}     
Let $T$ be a linearly repetitive tiling in $\RR^d$. Let $w$ be a subadditive, asymptotically $T$-invariant function. Then for every sequence $B_n$ with $B_n\in \boxen (r_n)$ and $r_n\to \infty $, the limit       
$$     
\lim_{n\to \infty} \frac{w(B_n)}{|B_n|}     
$$     
exists and is independent of the sequence.     
\end{coro}     
\noindent{\it Proof.} The strategy of the proof is simple. We will construct a subadditive function $F=F_w$ on $\Pbox(T)$ and show that the limit in question equals the limit of $F^+(r)$, whose existence is guaranteed by Theorem \ref{subadditive}.   
   
Define $F$ on $\Pbox(T)$ by      
$$     
F(P)= \sup\{ w(s(\Pdot)): \Pdot= T\cap s(\Pdot), \mbox{$\Pdot$ representative of $P$}\}.     
$$     
By properties (i) and (ii) of $w$, the function  $F$ is subadditive on $\Pbox(T)$. By construction of $F$, we have     
$$     
\frac{w(B)}{|B|}\leq F^+(r)     
$$     
for $B\in \boxen (r)$ with $r\geq r_w$ arbitrary. This immediately implies     
     
\begin{equation}\label{ugleins}      
\limsup_{n\to \infty} \frac{w(B_n)}{|B_n|}\leq \limsup_{r\to \infty} F^+(r).     
\end{equation}     
Moreover, by property (iii) of $w$, we have     
$$     
\frac{w(B)}{|B|}+ e_w(\omega(B))\geq F^-(r)     
$$     
for $B\in \boxen (r)$. This yields     
      
\begin{equation}\label{uglzwei}      
\liminf_{n\to\infty}  \frac{w(B_n)}{|B_n|}\geq \liminf_{r\to \infty} F^-(r).     
\end{equation}     
By $(\ref{ugleins})$, $(\ref{uglzwei})$, and Theorem \ref{subadditive},     
the statement of the theorem follows. \hfill \qedsymbol

\begin{coro}\label{delonesubadditive}     
Let $D$ be a linearly repetitive Delone set in $\RR^d$. Let $w$ be a subadditive, asymptotically $D$-invariant function. Then for every sequence $B_n$ with $B_n\in B(r_n)$ and $r_n\to \infty $, the limit       
$$   
\lim_{n\to \infty} \frac{w(B_n)}{|B_n|}   
$$     
exists and is independent of the sequence.     
\end{coro}   
\noindent{\it Proof.} This follows from the foregoing corollary applied to the a colored version of the Voronoi construction $V(D)$ (cf.~Section \ref{Preliminaries}). Here, each tile in $V(D)$ is colored by the unique   point of $D$ in its interior. To emphasize the coloring, we denote the  colored tiling by $V(D,C)$. The coloring implies that the function $w$ is asymptotically $V(D,C)$-invariant. Thus, the result follows from the foregoing corollary.  \hfill \qedsymbol    
   
\medskip   
   
Let us now discuss two classes of examples of the above theorems. They are given by tilings arising from primitive substitutions and tilings arising from Sturmian dynamical systems whose rotation number has bounded continued fraction expansion.   
   
We start by considering primitive substitutions. They give rise to linearly repetitive tilings \cite{DZ,DHS,sol1}. Thus, we immediately get the following result.   
   
\begin{coro}     
Let $S$ be a primitive substitution and let $T$ be a tiling associated to $\Pset(S)$ with vertex set $E$. Let $w$ be an asymptotically $E$-invariant subadditive function on boxes in $\RR^d$. Then, the limit $\lim_{n\to \infty}\frac{w(B_n)}{|B_n|}$ exists for every sequence $B_n$ of boxes with $B_n \in B(r_n)$ and $r_n \to \infty$, and it is independent of the sequence.      
\end{coro}     
This is essentially the subadditive ergodic theorem of \cite{GH}. The theorem of \cite{GH} is slightly more general in that the sequences $(B_n)$ considered there are only required to be cube-like van Hove sequences. On the other hand, the notion of subadditivity used there is  more restrictive than the notion used here. There, $w$ is required to satisfy a subadditivity condition on unions of quite general disjoint (up to their boundary) sets with the constant $c_w$ being zero. In fact, under these assumptions, one should be able to extend our theorem to hold for arbitrary cube-like van Hove sequences. However, our theorem is good enough to cover the desired applications.   
   
The other example is given by certain Sturmian dynamical systems; see Appendix B for some background. As shown in \cite{LP}, a Sturmian dynamical system is linearly repetitive if and only if its rotation number has bounded continued fraction expansion. Thus, we immediately obtain the following corollary of Theorem \ref{subadditive} which generalizes Theorem 2 of \cite{len} (cf.~\cite{dl2,l3}) as well).     
   
\begin{coro}     
Let an irrational $\alpha \in (0,1)$ with bounded continued fraction expansion be given. Let $\worte(\alpha)$ be the set of pattern classes of the Sturmian dynamical system with rotation number $\alpha$ {\rm (}cf.~\cite{dl2,len} for details{\rm )}. Then for every subadditive function $F$ on $\worte(\alpha)$, the limit $\lim_{n\to \infty}\frac{F(w_n)}{|w_n|}$     
exists for every sequence $(w_n)$ with $|w_n|$  going to infinity. Moreover, the limit is independent of the sequence.      
\end{coro}     
     
Of course, one could use Corollary \ref{additive} instead of Theorem     
\ref{subadditive} to obtain an additive ergodic theorem. However, this     
kind of result falls clearly short of the additive theorem of \cite{LP}, as it does not allow one to estimate the rate of convergence. This has been discussed in Remark \ref{Konvergenzrate} in Section \ref{Uniform}.   
   
We close this section by sketching a direct derivation of the unique   
ergodicity of dynamical systems associated to linearly repetitive   
tilings. In fact, the result uses only the compactness of the underlying   
space and an additive ergodic theorem. Thus, the proof applies verbatim to more general systems. The ``inner box'' technique given below applies to several contexts (cf. \cite{l3} for further discussion). It will be used in the next section as well.    
   
\begin{coro}   
Let $\Pset$ be linearly repetitive. Then the tiling dynamical system $(\Omega(\Pset), \RR^d)$ is uniquely ergodic. Here, $\RR^d$ acts on $\Omega(\Pset)$ in the canonical way via translation.   
\end{coro}   
\noindent{\it Proof.} We have to show that for any continuous $f$ on $\Omega(\Pset)$, the limits $\frac{1}{B} \int_B f(T - t)\,dt$ converge   
uniformly in $T$ for $\omega(B)$ going to infinity. The strategy is   
similar to the proof of Corollary \ref{tilingsubadditive} above. We will   
associate to $f$ additive functions $F_{\rm sup}$ and $F_{\rm inf}$   
on $\Pset$. They are defined as follows:   
$$   
F_{\rm sup} (P)=\sup\left\{\int_{s(\Pdot)} f(T-t)\,dt: T\cap s(\Pdot)=   
\Pdot\right\},   
$$   
$$   
F_{\rm inf} (P)=\inf\left\{\int_{s(\Pdot)} f(T-t)\,dt: T\cap s(\Pdot)=   
\Pdot\right\},   
$$   
where $\Pdot$ is an arbitrary representative of $P$ (it is not hard to   
check that these definitions are independent of the actual choice of   
$\Pdot$). Apparently,   
    
\begin{equation}\label{sub} F_{\rm sup} \left( \bigoplus_{j=1}^n P_j \right) \leq \sum_{j=1}^n F_{\rm sup}(P_j), \;\:\;F_{\rm inf} \left( \bigoplus_{j=1}^n  P_j \right) \geq \sum_{j=1}^n F_{\rm inf}(P_j).   
\end{equation}   
Moreover, the following is valid,   
   
\begin{equation} \label{littleo}   
|F_{\rm sup}(P)- F_{\rm inf}(P)| \leq o(\omega(P)),   
\end{equation}   
where the little $o$ function only depends on the continuity properties   
of $f$. To prove (\ref{littleo}) we use an ``inner box'' argument. Recall   
that $f$ is continuous and thus uniformly continuous since $\Omega(\Pset)$ is compact by Proposition \ref{compact}. Therefore, for each $\epsilon >0$, there exists $R$ such that $|f(T)-f(S)|\leq \epsilon$ whenever $T\cap B(0,R)=S\cap B(0,R)$. This implies that for all $t\in s(\Pdot)$ with $\dist(t,(s(\Pdot)^c)\geq R$, the difference of the integrands $|f(T - t)-f(S - t)|$ is smaller than $\epsilon $.  For large enough $\omega(P)$, the set of those $t$ agrees with the size of $P$ up to a boundary term. This proves (\ref{littleo}). By (\ref{sub}) and  (\ref{littleo}), the functions $F_{\rm sup}$ and $F_{\rm inf}$ are additive. Thus, the additive ergodic theorem implies the existence of the limits $\lim_{\omega(P)\to\infty} \frac{F_{\rm sup}(P)}{|P|}$ and $\lim_{\omega(P)\to\infty} \frac{F_{\rm inf}(P)}{|P|}$. By (\ref{littleo}), the limits are equal and the corollary follows. \hfill \qedsymbol

\section{Applications}\label{Applications}   
In this section we consider applications to random operators associated to tilings. In this context, there are two important quantities whose existence is established by a subadditivity argument, namely, the Lyapunov exponent in the one-dimensional case and the integrated density of states in arbitrary dimensions.   
   
The existence of the integrated density of states for Schr\"odinger-type   
operators associated to primitive substitutions is thoroughly discussed   
in  \cite{hof}. The discussion given there relies on abstract operator   
theory together with a subadditive ergodic theorem. Thus, it gives   
essentially the existence of the integrated density of states for   
Schr\"odinger-type operators associated to arbitrary linearly repetitive   
structures. In fact, the argument of \cite{hof} can be improved and   
strengthened in several respects \cite{l2,l3}. In particular, it turns   
out that the existence proof can actually be reduced to an additive   
ergodic theorem. This is interesting due to the existence of an error   
estimate in the additive ergodic theorem. This might have useful   
applications.    
   
Let us be more precise. For a tiling or pattern $M$, the space $l^2(M)$ is defined to be the space of all square summable sequences indexed by the elements of $M$. Let $A$ be a selfadjoint operator on a linearly repetitive tiling $T$ with matrix elements $A(x,y)$ for $x,y\in T$.  (Here, the tiling $T$ is called linearly repetitive if $\Pset(T)$ is linearly repetitive.) We will assume that $A$ satisfies the following finite range (FR) and invariance (I) properties: There exists some $R\geq 0$ with   
   
\begin{itemize}   
\item [(FR)] $A(x,y)$ vanishes for $\dist(x,y)\geq R$.    
\item[(I)] The value of $A(x,y)$ is completely determined by the   
pattern class $[\{t\in T: \dist(t,\{x,y\})\leq R\}]$.   
\end{itemize}   
In fact, the invariance condition implies that the operator $A$ can be defined on every tiling $T$ of the species $\Omega(T)$. To emphasize this, we will sometimes write $A(T)$ for the manifestation of $A$ on $l^2(T)$.   
   
For a box $B$ in $\RR^d$, the restriction  $A(T)|_B$ of $A$ to $B$ is the   
operator on $l^2(B\cap T)$ with matrix elements   
$$   
A(T)|_B(\tilde{x},\tilde{y})= A(x,y),\;\:\; \mbox{for $\tilde{x}=x\cap T$  and $\tilde{y}=y\cap T$ }.   
$$   
For a box $B$ in $\RR^d$ and $\lambda\in \RR$, define the function  $k^T_{\lambda}(B)$ by   
$$   
k^T_{\lambda}(B)=\frac{1}{|B|} \, \#\{\lambda_n: \lambda_n\leq   
\lambda, \;\lambda_n\;\mbox{eigenvalues of}\; A(T)|_B\},   
$$   
where the number of elements of a finite set $S$ is denoted by $\#S$. Then the following holds.   
   
\begin{theorem}\label{ids}   
The limit $\lim_{\omega(B)\to \infty} k^T_{\lambda}(B)$ exists and is independent of $T$. In fact, the convergence {\rm (}in $\omega(B)${\rm )} is uniform in $T$.   
\end{theorem}   
\noindent{\it Proof (sketch).} By Corollary \ref{tilingsubadditive} it is enough to show that the map $B\mapsto k^T_{\lambda}(B)$ is translation-invariant and additive. But this follows from the finite range condition together with the invariance condition. Details can be found in \cite{l2,l3}. \hfill \qedsymbol   
   
\medskip   
   
Theorem \ref{ids} generalizes the corresponding theorem of \cite{hof},   
where Penrose tilings are considered. Moreover, it only relies on an   
additive ergodic theorem, whereas \cite{hof} uses a subadditive theorem.   
   
Let us now turn to the study of the Lyapunov exponent. The sketch     
below follows the detailed discussion of the Sturmian case in \cite{dl2}. An admissible set $\Pset$ of abstract patterns in one dimension over a finite set of tiles can easily be identified with a set $\worte$ consisting of finite words over a finite alphabet $A \subset \RR$. The study of one-dimensional Schr\"odinger operators associated to $\worte$ can be based on the study of the so-called transfer matrices. For each $E\in \CC$, the transfer matrix $M(E)$ gives a map $M(E):\worte \longrightarrow {\rm SL}(2,\CC)$, defined by $M(E)(w)= T(E, w_n)\times \cdots \times T(E,w_1)$ for $w=w_1\ldots w_n$, where for $a\in \RR$ and $E\in \CC$, the matrix $T(E,a)$ is defined by     
     
\begin{equation}       
T( E,a)=\left( \begin{array}{cc} E- a & -1\\1 & 0 \end{array}        
\right).        
\end{equation}       
This map is antimultiplicative if the operation on $\worte$ is standard  concatenation of words. Since the standard norm $\|\cdot\|$ on ${\rm SL}(2,\CC)$ is submultiplicative, the function     
$$     
F:\worte\longrightarrow \RR, \;\,F(w)= \ln\|M(w)\|     
$$     
is subadditive. Thus, the results of Section \ref{Uniform} give the     
following theorem.     
   
\begin{theorem}     
Let $\worte$ and $F$ be as above. If $\worte$ is linearly repetitive, then the limit $\lim_{n\to \infty} \frac{F(w_n)}{|w_n|}$ exists for each sequence $(w_n)$ in $\worte$ with $|w_n|$ going to infinity and the limit does not depend on the sequence.      
\end{theorem}     
The limit in the theorem is called the Lyapunov exponent. It plays an     
important role in the study of one-dimensional Schr\"odinger operators.   
The theorem applies in particular to systems arising from primitive   
substitutions and to Sturmian dynamical systems whose rotation number   
has bounded continued fraction expansion.  This is due to the fact that   
these systems are linearly repetitve, as discussed in Section   
\ref{Specializing} (cf.~\cite{DZ,DHS,sol2} and \cite{LP} as well). Thus,   
the theorem generalizes the corresponding theorems of \cite{hof} and   
\cite{dl2, len}.   
   
Let us close this section by pointing out that there is a theory     
of lattice gas models for tilings arising from primitive substitutions     
\cite{GH}. This theory is built upon the subadditive theorem of     
\cite{GH}. Thus, it is very likely that considerable portions of it can     
be carried over to gas models on linearly repetitive tilings.      
\medskip     
     
\noindent        
{\it Acknowledgments.} D.~D.~was supported by the German Academic         
Exchange Service through Hochschulsonderprogramm III (Postdoktoranden)     
and  D.~L.~received financial support from Studienstiftung des Deutschen Volkes (Doktorandenstipendium), both of which are gratefully acknowledged.           
     
\begin{appendix}   
   
\section{Compactness of Linearly Repetitive Tiling Spaces}   
   
In this section we sketch a proof of Proposition \ref{compact}. It consists of two steps, namely, establishing a finiteness condition and performing a standard diagonalization procedure; compare \cite{rw,r}.   
   
Let $\Pset$ be a linearly repetitive admissible set of abstract patterns. We want to show that $\Omega(\Pset)$ is compact. Observe first that for every $r \ge 0$, there are only finitely many pattern classes $P \in \Pset$ with $r_{{\rm out}}(P) \le r$. To see this, consider any pattern class $Q$ such that $K(0,c_{{\rm LR}} \cdot r) \subset s(\Qdot)$ for some representative $\Qdot$ of $Q$. Delete from $\Qdot$ all the tiles which have empty intersection with $K(0,c_{{\rm LR}} \cdot r)$ and call the resulting pattern $\Qdot '$ and its pattern class $Q'$. It is clear that $\Qdot '$ has finite volume and that $Q'$ contains every abstract pattern $P \in \Pset$ with $r_{{\rm out}}(P) \le r$. This proves the assertion.   
   
Let us now consider a sequence $(T_n)_{n \in \NN}$ in $\Omega(\Pset)$. It suffices to prove that $(T_n)_{n \in \NN}$ has a convergent subsequence. To find this subsequence, we will inductively define sequences $(T_n^m)_{n \in \NN}$, $m \in \NN$, such that $(T_n^1)_{n \in \NN}$ is a subsequence of $(T_n)_{n \in \NN}$ and for $m_1 \ge m_2$, $(T_n^{m_1})_{n \in \NN}$ is a subsequence of $(T_n^{m_2})_{n \in \NN}$. This will be done in a way such that $(T_n^n)_{n \in \NN}$ converges. Choose any monotonically increasing sequence $r_m \rightarrow \infty$. Essentially, we will force the sequence $(T_n^m)_{n \in \NN}$ to converge on $K(0,r_m)$. It is then obvious from the definition of $d(\cdot,\cdot)$ that the diagonal sequence $(T_n^n)_{n \in \NN}$ will be $d$-Cauchy with obvious limit tiling.   
   
To define the refinement $(T_n^m)_{n \in \NN}$ of $(T_n^{m-1})_{n \in \NN}$ (think of $(T_n)_{n \in \NN}$ as $(T_n^0)_{n \in \NN}$), we will proceed in two steps. First, consider the pattern $\Pdot_n$ of tiles in $T_n^{m-1}$ having nonempty intersection with $K(0,r_m)$. The patterns $\Pdot_n$ have outer radius bounded by $r_m + 2 r_{{\rm max}}$ (with $r_{{\rm max}}$ from Definition 2.1) and hence, by the above observation, their abstract pattern classes $P_n$ belong to a finite set. Hence, one of them, say $P$, occurs infinitely often. Delete all the tilings from the sequence $(T_n^{m-1})_{n \in \NN}$ which have $P_n \not= P$. By the Selection Theorem \cite{GS,rw}, the remaining sequence has a subsequence such that the corresponding sets $\Pdot_{n_k}$ converge with respect to standard Hausdorff metric. Call this sequence $(T_n^m)_{n \in \NN}$. By the above remarks, it is easy to see that $(T_n^n)_{n \in \NN}$ is a convergent subsequence of $(T_n)_{n \in \NN}$.

\section{Strictly Ergodic Subshifts for Which the Uniform Subadditive Ergodic Theorem Fails}   
   
In this section we present one-dimensional examples which show that our main result fails if we only require strict ergodicity rather than linear repetitivity. We will consider a standard symbolic form of Sturmian tilings of the real line, that is, we study two-sided sequences over the alphabet $A=\{0,1\}$; see \cite{b,l} for background on Sturmian sequences.   
   
Let us first recall some standard notation. Given a finite alphabet $A$, we denote by $A^*$ the set of finite words over $A$ and by $A^\NN$ (resp., $A^\ZZ$) the set of one-sided (resp., two-sided) sequences over $A$, both called infinite words. Given a finite or infinite word $w$, we denote by ${\rm Sub}(w)$ the set of all finite subwords of $w$. Finally, given two finite words $v,w$, $\#_v(w)$ denotes the number of occurrences of $v$ in $w$.   
   
Fix some irrational $\alpha \in (0,1)$ and define the words $s_n$ over the alphabet $A$ by $s_{-1} = 1$, $s_0 = 0$, $s_1 = s_0^{a_1 - 1} s_{-1}$, and $s_n = s_{n-1}^{a_n} s_{n-2}$, $n \ge 2$, where the $a_n$ are the coefficients in the continued fraction expansion of $\alpha$. By definition, for $n \ge 2$, $s_{n-1}$ is a prefix of $s_n$. Therefore, the following (``right''-) limit exists in an obvious sense, $c_\alpha = \lim_{n \rightarrow \infty} s_n \in A^\NN$.        
   
Define the associated set of pattern classes $\worte(\alpha) \subset A^*$ by $\worte(\alpha) = {\rm Sub} (c_\alpha)$. The associated symbolic dynamical system $(\Omega(\alpha),T)$ is then given by $\Omega(\alpha) = \{ x \in A^\ZZ : {\rm Sub} (x) \subset \worte(\alpha)\}$ and $(Tx)_n = x_{n+1}$. $(\Omega(\alpha),T)$ is strictly ergodic for every irrational $\alpha$. It is linearly repetitive if and only if the sequence $(a_n)_{n \in \NN}$ is bounded.   
   
We will prove the following theorem.   
   
\begin{theorem}\label{counter}   
There exist $\alpha \in (0,1)$ irrational and a subadditive function $F$ on $\worte(\alpha)$ with the following property: There exist sequences $(w_n^k)_{n \in \NN}$ in $\worte(\alpha)$, $k = 1,2$, such that $|w_n^k| \rightarrow \infty$ as $n \rightarrow \infty$, $k = 1,2$, and   
$$   
\limsup_{n \rightarrow \infty} \frac{F(w_n^1)}{|w_n^1|} < \liminf_{n \rightarrow \infty} \frac{F(w_n^2)}{|w_n^2|}.   
$$   
In particular, the limit $\lim_{|w| \rightarrow \infty} \frac{F(w)}{|w|}$   
does not exist, that is, the uniform subadditive ergodic theorem does not hold for $\worte(\alpha)$.   
\end{theorem}   
   
The following properties of the words $s_n$ are well known and will be useful in the proof of Theorem \ref{counter}.   
    
\begin{prop}\label{wunderformel}     
{\rm (i)} For all $n\geq 2$, the word $s_n$ is a prefix of the word $s_{n-1} s_n$.\\   
{\rm (ii)} For every $n$, there is no nontrivial occurrence of $s_n$ in $s_n s_n$, that is, $s_n s_n = w_1 s_n w_2$ implies $w_1 = \varepsilon$ or $w_2 = \varepsilon$.   
\end{prop}      
   
We are now in position to give the   
   
\medskip   
{\noindent}{\it Proof of Theorem \ref{counter}.} Define the function $G$ on $\worte(\alpha)$ by   
$$   
G(w) = \sum_{n=1}^\infty \#_{s_{n-1} s_n} (w) (|s_{n-1}| + |s_n|).   
$$   
It is clear that all but finitely many of the terms are zero. Moreover, it is obvious that $G$ is superadditive. Thus, by Theorem 2 of \cite{len}, $\overline{G} = \lim_{n \rightarrow \infty} \frac{G(s_n)}{|s_n|}$ exists, but it is possibly infinite. Observe that $G(s_n)$ only depends on the numbers $a_1, \ldots, a_n$. Hence, by (using Proposition \ref{wunderformel})   
$$   
\frac{G(s_{n+1})}{|s_{n+1}|} \le \frac{a_{n+1} G(s_n) + G(s_{n-1}) + a_{n+1} \sum_{i=1}^{n-1}\left(|s_{i-1}| + |s_i|\right)}{a_{n+1} |s_n|},   
$$   
we see that we can force $\overline{G}$ to be finite if $a_n \rightarrow \infty$ sufficiently fast. But then we have   
$$   
\frac{G(s_{n-1} s_n)}{|s_{n-1} s_n|} = \frac{G(s_{n-1})}{|s_{n-1} s_n|} + \frac{G(s_n)}{|s_{n-1} s_n|} + \frac{|s_{n-1} s_n|}{|s_{n-1} s_n|} \ge \frac{G(s_n)}{|s_n| \left( 1 + \tfrac{|s_{n-1}|}{|s_n|} \right)} + 1,   
$$   
that is, $\frac{G(s_{n-1} s_n)}{|s_{n-1} s_n|}$ does not converge to $\overline{G}$. We can therefore conclude by setting $F = - G$, $w_n^1 = s_{n-1}s_n$, and $w_n^2 = s_n$. \hfill \qedsymbol   
    
\medskip   
   
\begin{remark}   
{\rm The proof of Theorem \ref{counter} actually provides an uncountable set of numbers $\alpha$ such that the uniform subadditive ergodic theorem fails for $\worte(\alpha)$. This set, however, has Lebesgue measure zero. By a more sophisticated argument (see \cite{l3}), one may prove this result for all $\alpha$'s obeying $\sum_{n=1}^\infty \frac{1}{a_n a_{n+1}} < \infty$. Since this set still has Lebesgue measure zero (cf.~methods in \cite{cfs}), it may be interesting to establish results for the Lebesgue-generic set of $\alpha$'s with intermediate $(a_n)$-behavior.}   
\end{remark}   
   
\end{appendix}


\begin{thebibliography}{10}     
   
\bibitem{b} J.~Berstel, Recent results in Sturmian words, in {\it Developments in Language Theory}, Eds.~J.~Dassow and A.~Salomaa, World Scientific, Singapore (1996), 13--24   
   
\bibitem{cfs} I.~P.~Cornfeld, S.~V.~Fomin and Ya.~G.~Sinai, {\it Ergodic Theory}, Springer-Verlag, New York, Heidelberg, Berlin (1982)   
     
\bibitem{dl2} D.~Damanik and D.~Lenz, Uniform spectral properties of      
one-dimensional quasicrystals, II. The Lyapunov exponent, preprint (math-ph/9905008, mp-arc/99-184), to appear in {\it Lett.~Math.~Phys.}         
     
\bibitem{DZ} D.~Damanik and D.~Zare, Palindrome complexity bounds for      
primitive substitution sequences, to appear in {\it Discrete Math.}     
     
\bibitem{DHS} F.~Durand, B.~Host and C.~Skau, Substitution dynamical systems, Bratteli diagrams and dimension groups, {\it Ergodic Theory Dynam.~Systems} {\bf 19} (1999), 953--993     
     
\bibitem{GH} C.~Geerse and A.~ Hof, Lattice gas models on self-similar     
aperiodic tilings, {\it Rev.~Math.~Phys.} {\bf 3} (1991), 163--221     
     
\bibitem{GS} B.~Gr\"unbaum and G.~C.~Shephard, {\it Tilings and Patterns},     
Freeman and Company, New York (1987)     
     
\bibitem{hof} A.~Hof, Some remarks on discrete aperiodic Schr\"odinger operators, {\it J.~Stat.~Phys.} {\bf 72} (1993), 1353--1374       
     
\bibitem{LP} J.~C.~Lagarias and P.~A.~B.~Pleasants, Repetitive Delone sets     
and quasicrystals, preprint (math.DS/9909033)     
     
\bibitem{len} D.~Lenz, Hierarchical structures in  Sturmian dynamical systems, preprint       
   
\bibitem{l2} D.~Lenz, in preparation   
   
\bibitem{l3} D.~Lenz, {\it Aperiodische Ordnung und gleichm\"assige spektrale Eigenschaften von Quasikristal\-len}, Ph.~D.~Thesis, in preparation   
   
\bibitem{l} M.~Lothaire, {\it Algebraic Combinatorics on Words}, in   
preparation    
   
\bibitem{LuP} W.~Lunnon and P.~A.~B.~Pleasants, Quasicrystallographic     
tilings, {\it J.~Math.~Pures Appl.} {\bf 66} (1987), 217--263      
     
\bibitem{que} M.~Queff\'{e}lec, {\it Substitution Dynamical Systems - Spectral Analysis}, Lecture Notes in Mathematics, Vol.~1284, Springer, Berlin, Heidelberg, New York (1987)     
     
\bibitem{rad} C.~Radin, Miles of tiles, in {\it Ergodic theory of     
$\ZZ^d$ actions}, Proceedings of the Warwick symposium, Warwick,     
UK, 1993-1994, Eds.~M.~Pollicott et al., Lond.~Math.~Soc.~Lect.~Note Ser.~228, Cambridge University Press, Cambridge (1996), 237--258     
     
\bibitem{rw} C.~Radin and M.~Wolff, Space tilings and local isomorphism, {\it Geometriae Dedicata} {\bf 42} (1992), 355--360   
   
\bibitem{r} D.~J.~Rudolph, Markov tilings of $\RR^n$ and representations of $\RR^n$ actions, {\it Contemp.~Math.} {\bf 94} (1989), 271--290   
   
\bibitem{sen} M.~Senechal, {\it Quasicrystals and Geometry}, Cambridge     
University Press, Cambridge (1995)     
     
\bibitem{sol1} B.~Solomyak, Dynamics of self-similar tilings, {\it Ergodic Theory Dynam.~Systems} {\bf 17} (1997), 695--738     
     
\bibitem{sol2} B.~Solomyak, Nonperiodicity implies unique composition for self-similar translationally finite tilings, {\it Discrete Comput.~Geom.} {\bf 20} (1998), 265--279     
     
\end{thebibliography}
\end{document}